\documentclass[12pt,a4paper]{article}
\usepackage{amsmath}
\usepackage{amssymb}
\usepackage{t1enc}
\usepackage[latin1]{inputenc}
\usepackage[english]{babel}
\usepackage{amsfonts}
\usepackage{latexsym}
\usepackage{bm}
\usepackage{setspace}
\usepackage[toc,page]{appendix}
\usepackage{graphicx}
\usepackage{enumerate}
\usepackage[numbers,sort&compress]{natbib}
\usepackage{tikz}
\usepackage[font=small]{caption}

\usepackage{theoremref}
\usepackage{lineno}

\usepackage{mathtools}

\newtheorem{theorem}{Theorem}

\newtheorem{rem}{Remark}

\newtheorem{construction}{Construction}
\newtheorem{problem}{Problem}

\allowdisplaybreaks

\begin{document}
\title{Solution to a problem on sequences of clique isolation numbers 
\medskip\medskip
}

\author{Peter Borg \\[2mm]
\normalsize Department of Mathematics \\
\normalsize Faculty of Science \\
\normalsize University of Malta\\
\normalsize Malta\\
\normalsize \texttt{peter.borg@um.edu.mt}
\and
Alexander Farrugia \\[2mm]
\normalsize Department of Mathematics \\
\normalsize Junior College \\
\normalsize University of Malta\\
\normalsize Malta\\
\normalsize \texttt{alex.farrugia@um.edu.mt}
}

\date{}
\maketitle

\begin{abstract}
For any positive integer $k$ and any $n$-vertex graph $G$, let $\iota(G,k)$ denote the $k$-clique isolation number of $G$, that is, the size of a smallest subset $D$ of the vertex set of $G$ such that the closed neighbourhood $N[D]$ of $D$ intersects the vertex sets of the $k$-cliques contained by $G$ (equivalently, $G-N[D]$ contains no $k$-clique). Thus, $\iota(G,1)$ is the domination number of $G$, and $\iota(G,2)$ is the vertex-edge domination number of $G$. We solve the following problem, posed by Fenech, Kaemawichanurat and the present first author [Discrete Math.~343 (2020), 111879]. Which integer sequences can be realized as $\iota(G,1), \dots, \iota(G,k)$ for some graph $G$? We show that for any integer sequence $a_1, \dots, a_k$, a graph $G$ with $(\iota(G,1), \dots, \iota(G,k)) = (a_1, \dots, a_k)$ exists if and only if $a_1 \geq \dots \geq a_k \geq 0$. We also explicitly show that each realization is attained by a connected graph unless $a_2 = 0$ and $a_1 \geq 2$.
\end{abstract}

\section{Introduction} \label{Introsection}
The graph isolation problem was introduced by Caro and Hansberg in \cite{CaHa15, CaHa17}. It is a natural generalization of the classical domination problem \cite{C, CH, HHS, HHS2, HL, HL2}. Similarly to the historical development of domination, isolation has rapidly become a very active field of investigation. Many of the known isolation results are referenced in \cite{Borgrsc, Borgisdom}. The seminal paper \cite{CaHa17} addressed this problem from various angles and raised a number of questions, some of which are answered in \cite{Borg, BFK, Borgrsc}. In \cite{BFK}, Fenech, Kaemawichanurat and the present first author solved a problem in \cite{CaHa17} concerning isolation of cliques. They also posed a problem \cite[Problem~3.4]{BFK} concerning sequences of clique isolation numbers. We solve this problem.

The set of positive integers is denoted by $\mathbb{N}$. For $n \in \mathbb{N}$, $[n]$ denotes the set $\{1, \dots, n\}$ (that is, $\{i \in \mathbb{N} \colon i \leq n\}$). We take $[0]$ as the empty set $\emptyset$. Arbitrary sets are taken to be finite. For a set $X$, the set of $2$-element subsets of $X$ is denoted by ${X \choose 2}$ (that is, ${X \choose 2} = \{ \{x,y \} \colon x,y \in X, x \neq y \}$). For standard terminology in graph theory, we refer the reader to \cite{West}. Most of the terminology used here is defined in \cite{Borg}. 

Every graph $G = (V(G), E(G))$ is taken to be \emph{simple}, that is, its vertex set $V(G)$ and edge set $E(G)$ satisfy $E(G) \subseteq {V(G) \choose 2}$. We may represent an edge $\{v,w\}$ by $vw$. We call $G$ an \emph{$n$-vertex graph} if $|V(G)| = n$. 
For $v \in V(G)$, $N_{G}(v)$ denotes the set of neighbours of $v$ in $G$, $N_{G}[v]$ denotes the closed neighbourhood $N_{G}(v) \cup \{ v \}$ of $v$, and $d_{G}(v)$ denotes the degree $|N_{G}(v)|$ of $v$. For $X \subseteq V(G)$, $N_G[X]$ denotes the closed neighbourhood $\bigcup_{v \in X} N_G[v]$ of $X$, $G[X]$ denotes $(X, E(G) \cap {X \choose 2})$ (the subgraph of $G$ induced by $X$), and $G - X$ denotes $G[V(G) \setminus X]$ (the graph obtained by deleting the vertices in $X$ from $G$). Where no confusion arises, the subscript $G$ may be omitted; for example, $N_G(v)$ may be abbreviated to $N(v)$.  

A \emph{component of $G$} is a maximal connected subgraph of $G$. Clearly, the components of $G$ are pairwise vertex-disjoint, and their union is $G$. 

Consider two graphs $G$ and $H$. If $G$ is a copy of $H$, then we write $G \simeq H$. If $H$ is a subgraph of $G$, then we say that \emph{$G$ contains $H$}. 

For $k \in \mathbb{N}$, the graph $([k], {[k] \choose 2})$ is denoted by $K_k$. A copy of $K_k$ is called a \emph{$k$-clique} or simply a \emph{clique}. A $k$-clique contained by $G$ is called a \emph{$k$-clique of $G$}. If $D \subseteq V(G)$ such that $N[D]$ intersects the vertex sets of the $k$-cliques of $G$, then $D$ is called a \emph{$K_k$-isolating set of $G$} or a \emph{$k$-clique isolating set of $G$}. Note that $D$ is a $K_k$-isolating set of $G$ if and only if $G-N[D]$ contains no $k$-clique. Let $\iota(G, k)$ denote the size of a smallest $K_k$-isolating set of $G$. 

If $D \subseteq V(G) = N[D]$, then $D$ is called a \emph{dominating set of $G$}. Clearly, $D$ is a dominating set of $G$ if and only if $D$ is a $K_1$-isolating set of $G$. Thus, $\gamma(G) = \iota(G, 1)$. The size of a smallest dominating set of $G$ is called the \emph{domination number of $G$} and is denoted by $\gamma(G)$. Assume that $G$ is a connected $n$-vertex graph. One of the earliest domination results, due to Ore \cite{Ore}, is that $\gamma(G) \leq n/2$ if $G \not\simeq K_1$ (see \cite{HHS}). Deleting the closed neighbourhood of a dominating set yields the graph with no vertices. Deleting the closed neighbourhood of a $K_2$-isolating set yields a graph with no edges. In the literature, such a set is also called a \emph{vertex-edge dominating set} \cite{BCHH,Lewis,LHHF,Peters,TJ,Z}. Caro and Hansberg~\cite{CaHa17} proved that $\iota(G, K_2) \leq n/3$ unless $G \simeq K_2$ or $G$ is a $5$-cycle. This was independently proved by \.{Z}yli\'{n}ski \cite{Z} and solved a problem in \cite{BCHH}. Solving a problem posed in \cite{CaHa17} and generalizing both the bounds above, Fenech, Kaemawichanurat and the present first author~\cite{BFK} proved that $\iota(G, k) \leq n/(k+1)$ unless $G \simeq K_k$ or $k = 2$ and $G$ is a $5$-cycle (Ore's result is the case $k = 1$). It is worth pointing out that in \cite{BFK2}, they also showed that if $m = |E(G)|$ and $G \not\simeq K_k$, then $\iota(G, k) \leq (m+1)/({k \choose 2} + 2)$. Both bounds on $\iota(G,k)$ were shown to be sharp in the respective papers. As indicated above, many other isolation results followed (see \cite{Borgrsc, Borgisdom}).

A tuple whose entries are integers will be called an \emph{integer tuple}. We say that an integer tuple $(a_1, \dots, a_k)$ is \emph{clique isolation realizable}, or \emph{ci-realizable} in brief, if there exists a graph $G$ such that $(\iota(G,1), \dots, \iota(G,k)) = (a_1, \dots, a_k)$, in which case we also say that $(a_1, \dots, a_k)$ is \emph{ci-realized by $G$}. 

As mentioned above, we solve the following problem.

\begin{problem}[{\cite[Problem~3.4]{BFK}}] \label{problem} Which integer tuples are ci-realizable?
\end{problem}

\noindent We show that the ci-realizable integer tuples are the ones that are decreasing (monotonic non-increasing) and obviously have non-negative entries. However, the result we present goes beyond this. For problems concerning graph parameters, it is often desirable and important to consider connected graphs. Motivated by this, we particularly show that if ${\bf a} = (a_1, \dots, a_k)$ is ci-realizable, then ${\bf a}$ is ci-realized by a connected graph unless $a_1 \geq 2$ and $a_2 = 0$. 

\section{The result}

We now start working towards presenting and proving our result.

\begin{rem} \label{rem} \emph{Consider the exceptional case where ${\bf a} = (a_1, \dots, a_k)$ is ci-realizable, $k \geq 2$ and $a_2 = 0$, and let $G$ be a graph such that $\bf a$ is ci-realized by $G$. Since $a_2 = 0$, $G$ contains no copy of $K_2$, so $E(G) = \emptyset$. Thus, since $\iota(G,K_1) = a_1$, we have $|V(G)| = a_1$. Therefore, $G$ is a copy of the empty graph $([a_1], \emptyset)$, and hence, if $a_1 \geq 2$, then $G$ is disconnected.}
\end{rem}
\begin{construction} \label{const} \emph{Consider any integer tuple ${\bf a} = (a_1, \dots, a_k)$ such that $a_1 \geq \dots \geq a_k \geq 1$. Let $a_{k+1} = 0$. Let $B_1, \dots, B_{a_1}$ be pairwise vertex-disjoint cliques such that for each $i \in [k]$, $|V(B_j)| = i$ for each $j \in [a_i] \setminus [a_{i+1}]$. Note that $|V(B_1)| = k$ as $a_1 \geq \dots \geq a_k$. For each $j \in [a_1]$, let $w_{j,1}, \dots, w_{j,n_i}$ be the distinct vertices of $B_j$. Let $Q$ be a set such that $|Q| = a_1$ and $Q \cap V(B_j) = \emptyset$ for each $j \in [a_1]$, and let $v_1, \dots, v_{a_1}$ be the distinct elements of $Q$. Let $T$ be the graph with $V(T) = Q$ and $E(T) = \{v_{1}v_{j} \colon j \in [a_1] \setminus \{1\}\}$. Let $X = \{v_jw_{j,1} \colon j \in [a_1]\}$ and $Y = \{v_1w_{1,h} \colon h \in [k-1]\}$. Let $B_{\bf a}$ and $B_{\bf a}^{\rm c}$ be the graphs with 
\begin{gather} V(B_{\bf a}) = \bigcup_{j = 1}^{a_1} V(B_j), \quad E(B_{\bf a}) = \bigcup_{j = 1}^{a_1} E(B_j), \nonumber \\
V(B_{\bf a}^{\rm c}) = V(T) \cup V(B_{\bf a}), \quad E(B_{\bf a}^{\rm c}) = E(T) \cup X \cup Y \cup E(B_{\bf a}). \nonumber
\end{gather}
Clearly, $B_1, \dots, B_{a_1}$ are the components of $B_{\bf a}$, and $B_{\bf a}^{\rm c}$ is a connected graph.}
\end{construction}

\begin{figure}
    \centering
    \includegraphics[width=0.47\linewidth]{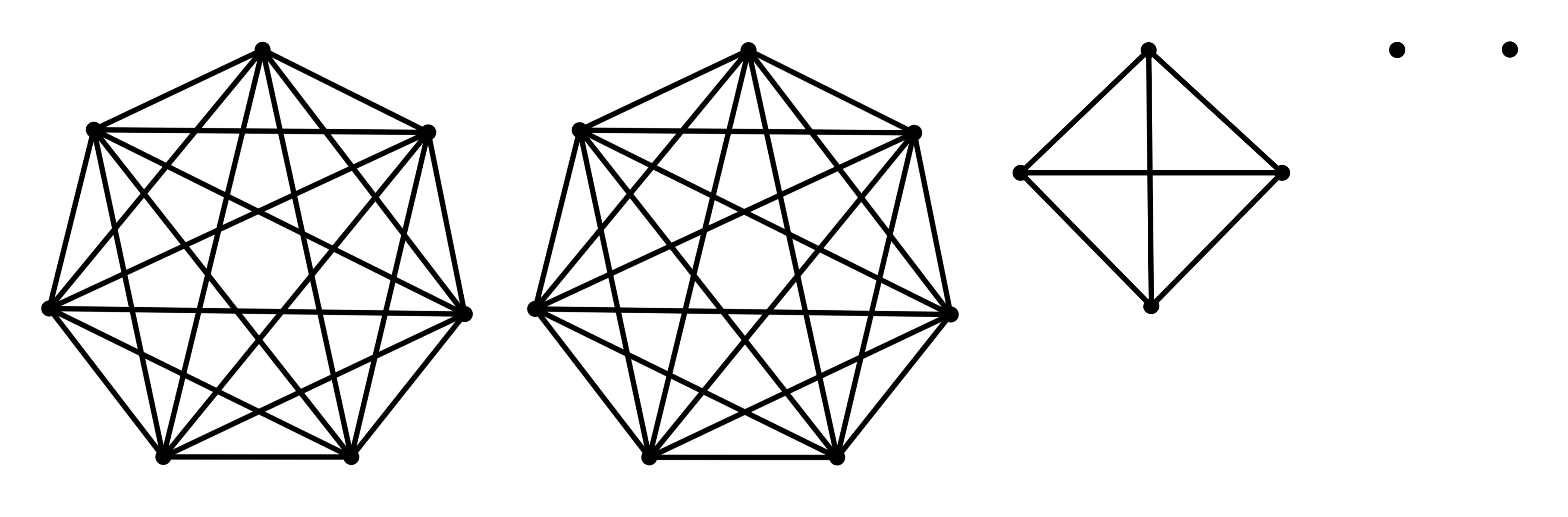} \quad \quad \includegraphics[width=0.47\linewidth]{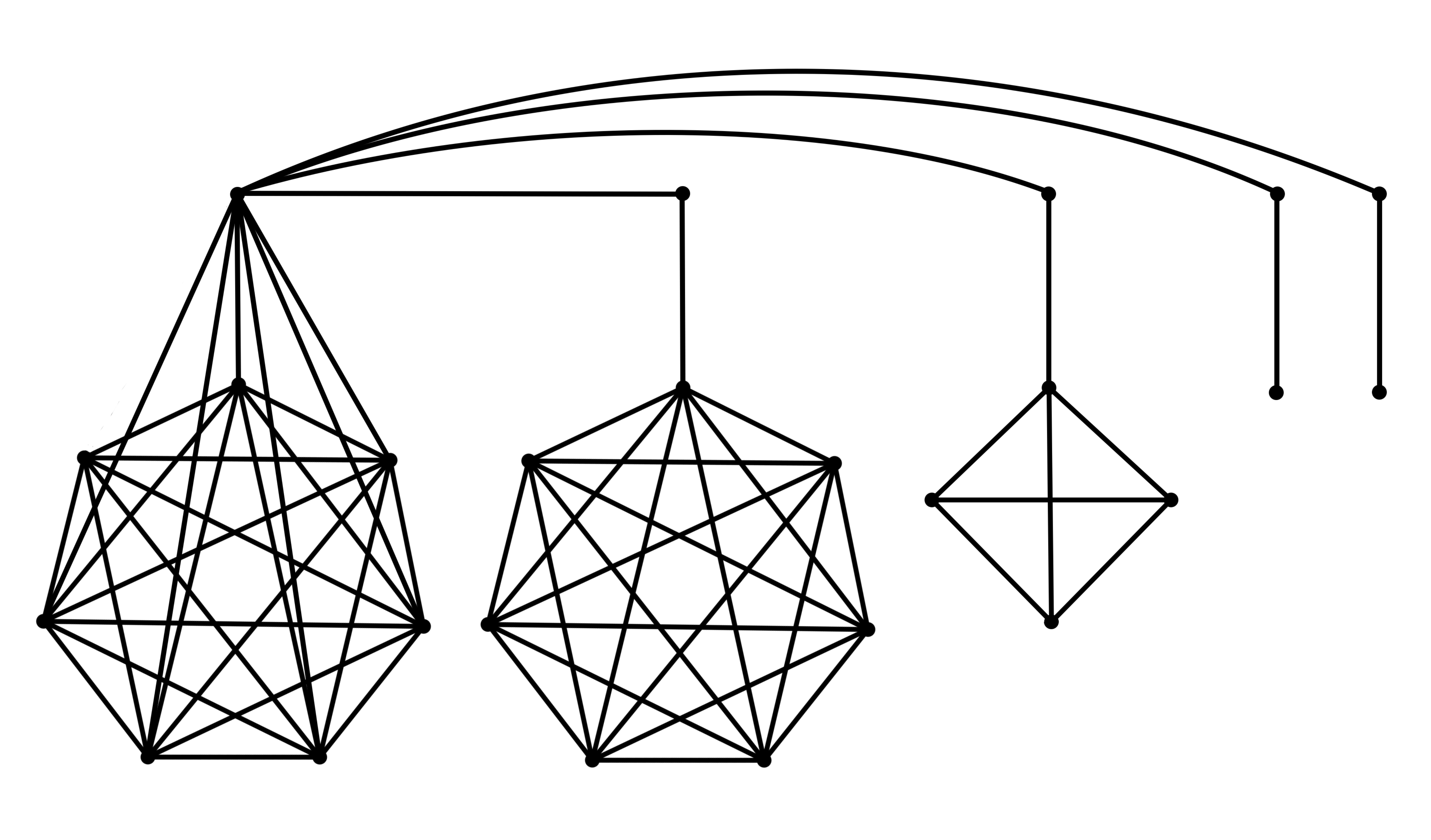}
    \caption{The graphs $B_{\bf a}$ (left) and $B_{\bf a}^{\rm c}$ (right) with ${{\bf a} = (5,3,3,3,2,2,2,0, \dots, 0)}$.}
    \label{fig:placeholder}
\end{figure}

Note that Construction~\ref{const} assumes that the entries of the tuple are all positive. Now consider any integer tuple ${\bf a} = (a_1, \dots, a_k)$ with $a_1 \geq \dots \geq a_k \geq 0$. If $a_1 = 0$, then let $B_{\bf a} = B_{\bf a}^{\rm c} = (\emptyset, \emptyset)$. In view of Remark~\ref{rem}, if $k \geq 2$ and $a_2 = 0$, then let $B_{\bf a} = B_{\bf a}^{\rm c} = ([a_1], \emptyset)$. If $a_2 \neq 0$, then let $\ell_{\bf a} = \max\{i \in [k] \colon a_i \neq 0\}$, let ${\bf a}' = (a_1, \dots, a_{\ell_{\bf a}})$, and let $B_{\bf a} = B_{{\bf a}'}$ and $B_{\bf a}^{\rm c} = B_{{\bf a}'}^{\rm c}$ (see Figure~\ref{fig:placeholder}).

\begin{theorem} \label{mainresult} An integer tuple ${\bf a} = (a_1, \dots, a_k)$ is ci-realizable if and only if 
\begin{equation} a_1 \geq \dots \geq a_k \geq 0. \label{dec}
\end{equation}
Moreover, if (\ref{dec}) holds, then the following statements hold: \\
(i) ${\bf a}$ is ci-realized by both $B_{\bf a}$ and $B_{\bf a}^{\rm c}$.\\
(ii) ${\bf a}$ is ci-realized only by disconnected graphs if and only if $k \geq 2$, $a_1 \geq 2$, $a_2 = 0$ and ${\bf a}$ is ci-realized only by empty $a_1$-vertex graphs. 
\end{theorem}
\textbf{Proof.} Suppose that ${\bf a}$ is ci-realizable. Then, there exists a graph $G$ such that $a_i = \iota(G,i) \geq 0$ for each $i \in [k]$. Let $D$ be a smallest $K_i$-isolating set of $G$. Suppose that $G$ contains a $K_{i+1}$-clique $H$. Let $W \subset V(H)$ with $|W| = i$. Then, $G[W]$ is a $k$-clique of $G$, so $\emptyset \neq N[D] \cap W \subseteq N[D] \cap V(H)$. Therefore, $D$ is a $K_{i+1}$-isolating set of $G$, and hence we have $\iota(G,i+1) \leq |D| = \iota(G,i)$. This yields (\ref{dec}).

We now prove the converse. Thus, suppose that (\ref{dec}) holds. We prove that (i) holds. If $a_1 = 0$, then $a_1 = \dots = a_k = 0$, $B_{\bf a} = B_{\bf a}^{\rm c} = (\emptyset, \emptyset)$, and hence (i) holds trivially. Suppose $a_1 \geq 1$. Let $G \in \{B_{\bf a}, B_{\bf a}^{\rm c}\}$. Let $W_i = \{w_{1,1}, \dots, w_{a_i,1}\}$ for each $i \in [k]$.

Suppose $k = 1$. Then, $V(B_j) = \{w_{j,1}\}$ for each $j \in [a_1]$. Clearly, $W_1$ is a $K_1$-isolating set of $G$, so $\iota(G,K_1) \leq a_1$. For each $j \in [a_1]$, $N[V(B_j)] \subseteq V(B_j) \cup \{v_j\}$. Thus, if $D$ is a $K_1$-isolating set of $G$, then $D \cap (V(B_j) \cup \{v_j\}) \neq \emptyset$ for each $j \in [a_1]$, and hence $|D| \geq a_1$. Together with $\iota(G,K_1) \leq a_1$, this gives us $\iota(G,K_1) = a_1$, as required. 

Now suppose $k \geq 2$. If $a_2 = 0$, then $a_2 = \dots = a_k = 0$, $G = ([a_1], \emptyset)$, $[a_1]$ is the only $K_1$-isolating set of $G$, and hence (i) holds. Suppose $a_2 \geq 1$.  

Suppose $a_k \geq 1$. By the argument for the case $k = 1$, $\iota(G,1) = a_1$. Now consider any $i \in [k]$ with $i \geq 2$, and let $A = [a_1] \setminus [a_i]$. If $G = B_{\bf a}$, then let $W = W_i$ and $Z = \emptyset$. If $G = B_{\bf a}^{\rm c}$, then let $W = (W_i \setminus \{w_{1,1}\}) \cup \{v_1\}$ and $Z = \{w_{1, k}\}$. We have $|V(B_j)| \geq i$ for each $j \in [a_i]$, and $|V(B_j)| < i$ for each $j \in A$. Let $H = G - N[W]$. Then, $V(H) = Z \cup \bigcup_{j \in A} V(B_j)$ and $E(H) = \bigcup_{j \in A} E(B_j)$ (recall that $V(B_1) = \{w_{1,1}, \dots, w_{1,k}\}$, so $V(B_1) \setminus N_{B_{\bf a}^{\rm c}}[v_1] = \{w_{1,k}\}$). Thus, since $|Z| < 2 \leq i$ and $|V(B_j)| < i$ for each $j \in A$, $W$ is a $K_i$-isolating set of $G$, and hence $\iota(G,K_i) \leq a_i$. For each $j \in [a_1]$, $N[V(B_j)] \subseteq V(B_j) \cup \{v_j\}$. Thus, if $D$ is a $K_i$-isolating set of $G$, then $D \cap (V(B_j) \cup \{v_j\}) \neq \emptyset$ for each $j \in [a_i]$, and hence $|D| \geq a_i$. Together with $\iota(G,K_i) \leq a_i$, this gives us $\iota(G,K_i) = a_i$, as required.

Now suppose $a_k = 0$. Recall that $a_2 \geq 1$, so $B_{\bf a} = B_{{\bf a}'}$ and $B_{\bf a}^{\rm c} = B_{{\bf a}'}^{\rm c}$. Let $h = \ell_{\bf a}$. We have $a_1 \geq \dots \geq a_h > a_{h+1} = \dots = a_k = 0$. By the argument for the case $a_k \geq 1$, $\iota(G,K_i) = a_i$ for each $i \in [h]$. If $G = B_{\bf a}$, then for each $i \in [k] \setminus [h]$, we trivially have $\iota(G,K_i) = 0$, so again $\iota(G,K_i) = a_i$, as required. Suppose $G = B_{\bf a}^{\rm c}$. Let $S \subseteq V(G)$ with $|S| = h+1$. Since $a_2 \geq 1$, we have $h \geq 2$, so $|S| \geq 3$. If $S \subseteq V(T)$, then $G[S] \not\simeq K_{h+1}$. We have $h = |V(B_1)| \geq \dots \geq |V(B_h)| \geq 1$, and since $w_{1,h} \notin N(v_1)$, $G[V(B_1) \cup \{v_1\}] \not\simeq K_{h+1}$. Thus, if $S \subseteq V(B_j) \cup \{v_j\}$ for some $j \in [h]$, then $G[S] \not\simeq K_{h+1}$. If $S \nsubseteq V(T)$ and $S \nsubseteq V(B_j) \cup \{v_j\}$ for each $j \in [h]$, then there exists some $j' \in [h]$ such that $x_1 \in V(B_j')$ and $x_2 \in V(G) \setminus (V(B_{j'}) \cup \{v_{j'}\})$ for some $x_1, x_2 \in S$, and hence, since $x_1x_2 \notin E(G)$, $G[S] \not\simeq K_{h+1}$. Therefore, for each $i \in [k] \setminus [h]$, $\iota(G,K_i) = 0 = a_i$, as required. 

It remains to prove (ii). Suppose that (\ref{dec}) holds. If $k \geq 2$, $a_1 \geq 2$ and $a_2 = 0$, then by Remark~\ref{rem}, ${\bf a}$ is ci-realized only by disconnected graphs. Suppose that ${\bf a}$ is ci-realized only by disconnected graphs. By (i), ${\bf a}$ is ci-realized by $B_{\bf a}^{\rm c}$, so $B_{\bf a}^{\rm c}$ is disconnected, and hence $k \geq 2$, $a_1 \geq 2$ and $a_2 = 0$. By Remark~\ref{rem}, ${\bf a}$ is ci-realized only by empty $a_1$-vertex graphs.~\hfill{$\Box$}

\end{document}